\documentclass[12pt]{article}
\usepackage{amsfonts}
\usepackage{amssymb}
\usepackage{amsmath}
\usepackage{fancyhdr}
\usepackage{eufrak}
\usepackage{latexsym}
\usepackage{amsbsy}
\usepackage[all]{xy}
\usepackage[latin1]{inputenc}
\usepackage[T1]{fontenc}        
\usepackage{amsmath}

\newtheorem{defn}{D\'efinition}[section]
\newtheorem{thm}[defn]{Th\'eor\`eme}
\newtheorem{prop}[defn]{Proposition}
\newtheorem{lem}[defn]{Lemme}
\newtheorem{cor}[defn]{Corollaire}

\newcommand{\proof}{\vskip 2mm \noindent {\textsc{D\'emonstration : }}\rm}
\newcommand{\fin}{\hfill{\Large$\Box$}\\}

\newcommand{\al}{\alpha}

\newcommand{\Ga}{\Gamma}

\newcommand{\pp}{\forall \, \hat x \ \hat \mu \textrm{-p.p.} \ , \ }
\newcommand{\si}{\sigma}

\newcommand{\om}{\omega}

\newcommand{\C}{\mathbb {C}}
\newcommand{\R}{\mathbb {R}}
\newcommand{\N}{\mathbb {N}}
\newcommand{\Z}{\mathbb {Z}}

\newcommand{\Pj}{\mathbb {P}}

\newcommand{\diam}{{\rm diam \,}}
\newcommand{\vol}{{\rm vol \, }}
\newcommand{\Lip}{{\rm Lip \, }}
\newcommand{\Jac}{{\rm Jac \,  }}
\newcommand{\dist}{{\rm dist \, }}

\def\bignorm#1{\left\|\, #1\,\right\|}

\def\AA{{\cal A}}
\def\BB{{\cal B}}
\def\CC{{\cal C}}
\def\PP{{\cal P}}

\def\GG{{\cal G}}
\def\JJ{{\cal J}}

\def\II{{\cal I}}

\def\VV{{\cal V}}

\def\GG{{\cal G}}

\def\HH{{\cal H}}
\def\com{\ar@{}[rd]|{\circlearrowleft}}

\title {Dimension de la mesure d'\'equilibre d'applications m\'eromorphes}
\author{Tien-Cuong Dinh et Christophe Dupont}
\date{ \today }

\begin{document}

\maketitle

\begin{abstract}
Let $f$ be a dominating meromorphic self-map of a compact K\"ahler manifold.
Assume that the topological degree of $f$ is larger than the other dynamical degrees. 
We give estimates of 
the dimension of the equilibrium measure of $f$, which involve the
Lyapounov exponents.  
\end{abstract}

2000 Mathematics Subject Classification : 37C45, 37F10, 32H50.

Key Words : dimension theory, equilibrium measure, Lyapounov exponent.


\section{Introduction}\label{intro}


Soient $X$ une vari\'et\'e k\"ahl\'erienne compacte de dimension $k$ et $f : X \to X$ une application m\'eromorphe dominante. On note $d_t$ son 
degr\'e topologique, $\lambda_{k-1}$ son $(k-1)$-i\`eme degr\'e dynamique 
(cf section \ref{nota}), et on suppose que $d_t > \lambda_{k-1}$.  
Fixons aussi une forme de K\"ahler $\om$ sur $X$, normalis\'ee par $\int_X \om^k = 1$. Sous ces hypoth\`eses, Guedj a montr\'e dans \cite{G} que la suite de mesures : 
\begin{equation*}
 \mu_n := {1 \over d_t ^n} {f^n}^* \om^k 
\end{equation*}
converge vers une mesure de probabilit\'e invariante $\mu$ (cf aussi \cite{RS}, \cite{DS1} et \cite{DS2}). Ce r\'esultat avait \'et\'e \'etabli auparavant pour les applications holomorphes 
de degr\'e alg\'ebrique $d\geq 2$ sur l'espace projectif $\Pj^k$ (alors $d_t = d^k$ et $\lambda_{k-1} = d^{k-1}$) et pour les applications d'allure polynomiale. On consultera \`a ce sujet les livres de 
Berteloot-Mayer \cite{BM} (en dimension 1) et de Sibony \cite{S} (en dimension sup\'erieure). On trouvera aussi dans ces ouvrages une bibliographie plus d\'etaill\'ee.
\\

La mesure limite $\mu$ int\`egre les fonctions \emph{quasi-psh}, ce qui assure l'existence de ses exposants de Lyapounov. Nous les noterons en ordre croissant $\chi_1 \leq \cdots \leq \chi_k$. Si $\Sigma$ d\'esigne leur somme, on dispose des in\'egalit\'es (cf \cite{BD}, \cite{DS1}, \cite{G}) : 
\[\chi_1 \geq {1 \over 2} \log {d_t \over \lambda_{k-1}} \ \  \textrm{ et } \ \ 2\Sigma \geq  \log d_t. \]

D'une mani\`ere g\'en\'erale, les exposants de Lyapounov sont reli\'es \`a l'entropie et \`a la dimension de $\mu$, not\'ees respectivement $h(\mu)$ et $\dim_\HH(\mu)$ (cf section \ref{dimhaus}). Par exemple, si $f$ d\'esigne une fraction rationnelle de degr\'e $d$ sur $\Pj^1$, on a $h(\mu) = \log d = \dim_\HH(\mu) . \chi$, o\`u $\chi$ est l'unique exposant de $\mu$ (cf \cite{L}, \cite{M}). Une telle formule provient du caract\`ere conforme de $f$. Notons qu'il existe une relation semblable pour les diff\'eomorphismes des surfaces r\'eelles compactes \cite{Y}. 

En dimension complexe plus grande que $1$, les applications m\'eromorphes ne sont plus conformes. On conna\^it une minoration de la dimension lorsque $\mu$ est une masse de Monge-Amp\`ere \`a potentiel h\"old\'erien (par exemple si $f$ est holomorphe sur $\Pj^k$, \cite{B} §III.2, \cite{S} §1.7). Cette minoration fait appara\^itre l'exposant de H\"older. Notons que $\mu$ ne poss\`ede pas toujours cette propri\'et\'e de r\'egularit\'e. On dispose aussi d'une majoration faisant intervenir les exposants lorsque $f$ est un endomorphisme holomorphe de $\Pj^k$ d'origine polynomiale \cite{BdM}. \\

Le but de cet article est d'\'etablir l'encadrement suivant :
\begin{thm}\label{result}
Sous les hypoth\`eses pr\'ec\'edentes, la dimension de $\mu$ v\'erifie :
\begin{equation}\label{qqq}
  {\log d_t \over \chi_k} \leq \dim_\HH (\mu) \leq 2k - { 2\Sigma -  \log d_t \over \chi_k }  .
\end{equation}
\end{thm}

La d\'emonstration repose essentiellement sur l'existence et le contr\^ole 
des branches inverses des it\'er\'es de $f$ le long d'orbites n\'egatives typiques (cf la proposition \ref{lempr}). Nous donnons dans le cas m\'eromorphe un \'enonc\'e semblable au cas holomorphe, d\^u \`a Briend-Duval \cite{BD}. Le nouvel ingr\'edient est le contr\^ole de la vitesse d'approche des orbites n\'egatives vers l'ensemble $\JJ$ qui contient les ensembles d'ind\'etermination de $f$ et $f^{-1}$. L'int\'egrabilit\'e de la fonction $\log (\dist(.,\JJ))$ (cf lemme \ref{dfc}) joue ici un r\^ole crucial. Notre \'enonc\'e apporte aussi une pr\'ecision suppl\'ementaire : il permet de comparer les familles de branches inverses $(f^{-n})_{n \geq 0}$ d\'efinies au voisinage de $f^p(x)$ lorsque $p$ varie. Il stipule par exemple que le rayon de la boule sur laquelle elles sont d\'efinies d\'ecro\^it exponentiellement lentement lorsque $p$ devient grand.  \\

Esquissons la preuve de la majoration de (\ref{qqq}). Nous reprenons en partie la m\'ethode de Binder et DeMarco \cite{BdM}. Il s'agit d'introduire un bor\'elien $A$ de mesure positive sur lequel on dispose d'un contr\^ole uniforme des branches inverses. La proposition \ref{lempr} nous fournit des suites $(r_n)_n$ et $(\gamma_n)_n$ \`a d\'ecroissance lente telles que pour tout $x \in A$ et pour tout $n \geq 0$ :

 - la branche inverse $g_n$ de $f^n$ v\'erifiant $g_n(f^n(x)) = x$ existe sur  
$\BB_n := B_{f^n(x)}(r_n)$,

- l'ouvert $g_n(\BB_n)$ contient $\BB_n' = B_x (\gamma_n e^{-n \chi_k})$,

- le volume de $g_n(\BB_n)$ est major\'e par $ce^{-2n\Sigma}$, o\`u $c$ est une constante.

\noindent On a not\'e $B_x(r)$ la boule ouverte centr\'ee en $x$ et de rayon $r$. Les boules du type $\BB_n'$ permettent alors de construire des recouvrements de $A$ de diam\`etre arbitrairement petit, et l'estimation sur le volume donne une borne sur leur cardinal. On obtient ainsi la majoration indiqu\'ee. \\  

Venons-en \`a la minoration. Elle repose sur le r\'esultat classique suivant (cf \cite{Y}) :
\begin{thm} \label{young} 
Soit $\Lambda$ un bor\'elien de $X$ tel que $\mu(\Lambda) > 0$. Supposons que pour tout $x \in \Lambda$, il existe une suite d\'ecroissante de r\'eels positifs $(\rho_n(x))_{n \geq 0}$ v\'erifiant :
\begin{equation}\label{sz}
 \lim_{n \to +\infty} \rho_n(x) = 0 \ , \  \lim_{n \to +\infty} {\log \rho_{n+1}(x) \over \log \rho_n(x)} = 1 
\end{equation}
et
\begin{equation}\label{xz}
 a \leq   \liminf_{n \to +\infty}{ \log \mu  (  B_{x} (\rho_n(x)) )  \over  \log  \rho_n(x) } \leq   \limsup_{n \to +\infty}{ \log \mu  (  B_{x} (\rho_n(x)) )  \over  \log  \rho_n(x)} \leq b.     
\end{equation}
Alors la dimension de Hausdorff de $\Lambda$ v\'erifie $a \leq  \dim_\HH (\Lambda) \leq b$.
\end{thm}
Notons que seule la minoration sera utilis\'ee, sous la forme du corollaire :
\begin{cor}\label{bore}
Si il existe un bor\'elien $\Lambda$ de $X$ v\'erifiant les propri\'et\'es (\ref{sz}) et (\ref{xz}), alors la dimension de $\mu$ satisfait :
\[ \dim_\HH(\mu) \geq a. \]   
\end{cor}
En effet, si $Z$ est un bor\'elien de $X$ de mesure $1$, alors les points de $\Lambda \cap Z$ satisfont (\ref{sz}) et (\ref{xz}), ce qui entra\^ine $\dim_\HH (Z) \geq \dim_\HH (\Lambda \cap Z) \geq a$.\\

L'ensemble $\Lambda$ sera en fait de mesure $1$, et proviendra de la proposition sur les branches inverses. Nous prendrons pour $B_{x} (\rho_n(x))$ les boules $\BB_n'$ introduites plus haut : la d\'ecroissance lente de leur rayon nous assure que (\ref{sz}) est bien satisfaite. La minoration dans (\ref{qqq}) est alors cons\'equence du corollaire. Signalons que la borne $a = {\log d_t \over \chi_k}$ obtenue dans le th\'eor\`eme provient du fait que $\mu$ est de jacobien constant $d_t$ (cf section \ref{nota}).\\   

Dans le cas holomorphe, le th\'eor\`eme \ref{result} se reformule ainsi :
\begin{cor}
Soit $f$ un endomorphisme holomorphe de $\Pj^k$ de degr\'e alg\'ebrique 
$d\geq 2$. Alors la dimension de la mesure d'\'equilibre $\mu$ de $f$ v\'erifie : 
\[  {k\log d \over \chi_k} \leq \dim_\HH (\mu) \leq 2k - { 2\Sigma - k\log d \over \chi_k  } \leq 2k-2 + {\log d \over \chi_k}.   \] 
\end{cor}
La derni\`ere majoration provient de l'in\'egalit\'e $\chi_1 \geq {1 \over 2} \log d$ due \`a Briend-Duval \cite{BD}. On retrouve alors l'estimation de Binder-DeMarco \'etablie pour les endomorphismes $f$ d'origine polynomiale \cite{BdM}. Lorsque les exposants de Lyapounov sont minimaux,
i.e. \'egaux \`a ${1 \over 2} \log d$, les in\'egalit\'es ci-dessus deviennent 
des \'egalit\'es, et seuls les exemples de Latt\`es v\'erifient cette propri\'et\'e
(cf \cite{BDu}, \cite{D}). Notons que notre r\'esultat reste valable pour les applications \`a allure polynomiale dont la mesure d'\'equilibre int\`egre les fonctions \emph{psh} \cite{DS1}. On obtient \'egalement ces estimations dans un cadre \emph{r\'eel}. Par exemple, il est possible d'associer \`a des rev\^etements ramifi\'es sur des vari\'et\'es riemaniennes des ``mesures d'\'equilibre'' $\mu$ naturelles (cf \cite{DS1} §2). Elles sont de jacobien constant $d_t$, et le jacobien de $f$ est $\mu$-int\'egrable (cf section \ref{nota}). Si les exposants de $\mu$ sont strictement positifs (cette hypoth\`ese est toujours satisfaite en dimension $1$), on d\'emontre de m\^eme la proposition \ref{lempr} sur les branches inverses et on en d\'eduit un encadrement analogue sur la dimension de $\mu$. \\

L'article s'organise de la mani\`ere suivante. La deuxi\`eme section consiste en des g\'en\'eralit\'es : nous y pr\'ecisons les notations, ainsi que la d\'efinition des fonctions \`a variation lente. La section suivante est consacr\'ee \`a la proposition sur les branches inverses. La d\'emonstration du th\'eor\`eme \ref{result} occupe finalement le dernier paragraphe.\\  

\noindent \emph{Remerciements :} Nous remercions le referee pour sa lecture attentive. Ses remarques judicieuses nous ont aid\'e \`a am\'eliorer la r\'edaction de l'article.  


\section{G\'en\'eralit\'es}\label{gene}



\subsection{D\'efinitions et remarques}\label{nota}


Nous travaillons avec la m\'etrique induite par une forme de K\"ahler $\om$ fix\'ee. 
Pour simplifier les estimations, nous l'identifions parfois dans les cartes avec la m\'etrique
standard. Une fonction $u$ sur $X$ est \emph{quasi-plurisousharmonique} (\emph{quasi-psh}) si elle est int\'egrable, semi-continue sup\'erieurement et v\'erifie $i\partial \bar \partial u \geq - c \om$, o\`u $c \geq 0$. Autrement dit, localement, $u$ est la somme d'une fonction \emph{psh} et d'une fonction lisse.\\ 

Pour tout $x \in X$ et $r > 0$, on note $B_x(r)$ la boule ouverte centr\'ee en $x$ et de rayon $r$. Le diam\`etre et le volume d'une partie $A$ de $X$ sont respectivement not\'es $\diam(A)$ et $\vol(A)$. La distance d'un point $x \in X$ \`a une partie $A$ est not\'ee $\dist(x,A)$. On note aussi $\Lip(h)$ la constante de Lipschitz d'une application $h : U \to X$ d\'efinie sur un ouvert $U$ de $X$.  \\

Soient $X_1$ et $X_2$ deux vari\'et\'es complexes compactes, $Z:=X_1 \times X_2$ et $\pi_1 , \pi_2$ les projections de $Z$ sur ses facteurs. Par d\'efinition, une application  m\'eromorphe multivalu\'ee $g : X_1 \to X_2$ est la donn\'ee de sous-ensembles analytiques $\Gamma$ de $Z$ et $\AA$ de $X_1$ tels que $\pi_1^{-1}(\AA)$ ne contient aucune composante de $\Gamma$ et tels que la restriction $\pi_1 : \Gamma \setminus \pi_1^{-1}(\AA) \to X_1 \setminus \AA$  est un rev\^etement non ramifi\'e de degr\'e $d$. En particulier, $\AA$ contient l'ensemble d'ind\'etermination de $g$ (\emph{i.e.} les points $x \in X_1$ tels que $\pi_1^{-1}(x) \cap \Gamma$ est infini) et les singularit\'es de $\Gamma$ sont contenues dans $\pi_1^{-1}(\AA)$. Au voisinage de tout point $z \notin \AA$, l'application $g$ est donc une collection de $d$ applications holomorphes $g_1,\ldots,g_d$ et on pose : \[  \bignorm {d g (z) }  :=  \max_{i \in \{ 1,\ldots,d \}} \bignorm {d g_i (z)} \ ,  \ \bignorm {d^2 g (z) }  :=  \max_{i \in \{ 1,\ldots,d \}}  \bignorm {d^2 g_i (z)}.    \] 
Pour $A\subset X_1$ et $B\subset X_2$, posons 
$g(A):=\pi_2(\pi_1^{-1}(A)\cap\Gamma)$ et $g^{-1}(B):=\pi_1(\pi_2^{-1}(B)\cap\Gamma)$.
Nous utiliserons le r\'esultat suivant :

\begin{lem}\label{cauchy}
Soient $g : X_1 \to X_2$, $g_i$ et $\AA$ comme ci-dessus. 
Il existe des constantes $\tau>0$ et $p \in \N^*$ telles que pour tout $x \in X_1 \setminus \AA$, on a 
\[  \bignorm { dg(x) } + \bignorm { d^2 g(x) } \leq \tau \, \dist(x, \AA)^{-p} \]
et
\[ \dist(g_i(x),g(\AA))\leq \tau \, \dist(x,\AA)^{1/p} .\] 
\end{lem}

\proof
V\'erifions que l'on peut supposer que $\Gamma$ est lisse. Posons \`a cet effet $\AA' := \pi_1 ^{-1} (\AA)$. D'apr\`es le th\'eor\`eme de d\'esingularisation d'Hironaka \cite{H}, il existe une vari\'et\'e complexe compacte $\widetilde Z$, une sous-vari\'et\'e lisse $\widetilde \Gamma$ de $\widetilde Z$ et une application holomorphe $\pi : \widetilde Z \rightarrow Z$ telles que $\pi$ soit d'une part une bijection
entre $\widetilde Z \setminus \pi^{-1} (\AA')$ et $Z\setminus \AA'$  et d'autre part une bijection
entre $\widetilde \Gamma \setminus \pi^{-1}(\AA')$ et $\Gamma \setminus \AA'$. Puisque l'application $\pi_1 \circ \pi : \widetilde \Gamma \to X_1$ est holomorphe, son inverse $\widetilde g : X_1 \rightarrow \widetilde \Gamma$ est une application m\'eromorphe multivalu\'ee. Ainsi, $g$ est la compos\'ee de $\widetilde g$ par l'application holomorphe $\pi_2  \circ \pi$ dont la diff\'erentielle est uniform\'ement born\'ee. Donc, quitte \`a remplacer $g$ par $\widetilde g$, on peut supposer que $\Gamma$ est lisse. 

Soient $(a,b) \in \Gamma$ et $M , N$ des cartes de $X_1 , X_2$ contenant respectivement $a$ et $b$. On les munit de coordonn\'ees 
$(x_1,\ldots,x_k)$ et $(y_1,\ldots,y_l)$ telles que $a=0$ et $b=0$. Soit $U \Subset M \times N$ un voisinage de $(0,0)$. Fixons $(x,y) \in \Gamma \cap U$ et notons $\BB \subset M \setminus \AA$ la boule centr\'ee en $x$ de rayon maximal telle que la composante $\Gamma_{(x,y)}$ de $\pi_1^{-1}(\BB) \cap \Gamma$ contenant $(x,y)$ soit incluse dans $M \times N$. Soit $g^*$ l'application univalente d\'efinie sur $\BB$ de graphe $\Gamma_{(x,y)}$. 

Montrons la premi\`ere in\'egalit\'e. Nous allons estimer les d\'eriv\'ees de $g^*$ en $x$ en utilisant la formule de Cauchy. Puisque $g^*$ est born\'ee (elle prend ses valeurs dans $N$), il suffit de montrer que $\BB$ contient une boule centr\'ee en $x$ et de rayon $\tau'' \dist(x,\AA)^q$, o\`u $\tau'' >0$ et $q \in \N^*$. Observons que sur la vari\'et\'e \emph{lisse} $\Gamma \cap M \times N$ de dimension $k$, les formes holomorphes $d y_j$ s'\'ecrivent de fa\c{c}on unique comme combinaisons lin\'eaires des formes $dx_i$ \`a coefficients m\'eromorphes.  Ces coefficients sont en fait holomorphes hors de $\AA'$ car $\Gamma$ y est localement une union de graphes. Ainsi, quitte \`a restreindre $M$ et $N$, on a $\bignorm { dg^*(z) } \leq \tau' \, \dist(z,\AA)^{-q}$ pour tout $z \in \BB$, o\`u $\tau' >0$ et $q \in \N^*$. On en d\'eduit, d'apr\`es l'in\'egalit\'e des accroissements finis, qu'il existe une boule centr\'ee en $x$ et de rayon $\tau'' \, \dist(x,\AA)^q$ contenue dans $\BB$. On utilise ici le fait que $\Gamma \cap M \times N$ est une sous-vari\'et\'e de $M \times N$ et le fait que $\pi_1 : \Gamma \setminus \AA' \to X_1 \setminus \AA$  est un rev\^etement non ramifi\'e. La premi\`ere in\'egalit\'e en d\'ecoule. 

Passons \`a la deuxi\`eme in\'egalit\'e. En posant $y := g^*(x)$, il suffit de montrer 
$\dist (y,g(\AA))\lesssim \dist(x,\AA)^{1/p}$. Soient $u$ et $v$ des
fonctions holomorphes vectorielles sur $M$ et $N$ telles que $\AA=\{u=0\}$ dans $M$ et 
$g(\AA) = \{v=0\}$ dans $N$. On dispose des estimations suivantes, avec 
$q \in \N^*$ :
\begin{equation}\label{loj1}
 \dist(x,\AA) ^q \lesssim \| u(x) \| \lesssim  \dist(x,\AA),
\end{equation}
\begin{equation}\label{loj2}
 \dist(y,g(\AA)) ^{q} \lesssim \| v(y) \| \lesssim  \dist(y,g(\AA)).
\end{equation}
Les in\'egalit\'es de gauche sont celles de Lojasiewicz (\cite{Lo}, IV.7.2) et celles de droite proviennent de l'in\'egalit\'e des accroissements finis. 
Observons que l'ensemble analytique  
$\AA_1:=\{z\in\Gamma\cap M\times N,\  u\circ \pi_1(z)=0\}$ 
est contenu dans $\AA_2:=\{z\in\Gamma\cap M\times N,\ v\circ\pi_2(z)=0\}$. 
La vari\'et\'e $\Gamma\cap M\times N$ \'etant lisse, on obtient les in\'egalit\'es suivantes de la m\^eme mani\`ere que (\ref{loj1}) et (\ref{loj2}) :
\begin{eqnarray*}  
\lefteqn{\|v(y)\|=\|v\circ\pi_2(x,y)\|\lesssim {\rm dist}_\Gamma ((x,y),\AA_2)} \\
&&
\leq{\rm dist}_\Gamma((x,y),\AA_1)\lesssim \|u\circ\pi_1(x,y)\|^{1/q'} = \|u(x)\|^{1/q'}. 
\end{eqnarray*}
Cela termine la d\'emonstration \fin

Dans toute la suite, $f : X \to X$ d\'esigne une application m\'eromorphe \emph{dominante}. Cela signifie que son jacobien (d\'efini plus bas) n'est identiquement nul sur aucun ouvert de $X$. Notons que l'application $f^{-1} : X \to X$ est m\'eromorphe et multivalu\'ee. Soient $d_t$ le {\it degr\'e topologique} de $f$, $\CC$ son ensemble critique et $\II$ son ensemble d'ind\'etermination. 
On pose $\JJ':= f(\CC\cup\II)$ et 
$\JJ := \JJ'\cup f^{-1}(\JJ')$. 
Dans la section  \ref{brinv}, nous appliquerons le lemme \ref{cauchy} \`a $g = f$ et $g =f^{-1}$ avec $\AA := \II$ ou $\JJ'$. \\
 
Pour tout $l \in \{1,\ldots,k\}$, on d\'efinit la forme $f^* \om^l$ comme l'extension triviale de $(f_{\vert X \setminus \II})^* \om^l$ \`a travers $\II$. Avec cette convention, $\Jac (f)$ est la fonction positive sur $X$ v\'erifiant 
$f^* \om^k = \Jac (f). \om^k$. Les {\it degr\'es dynamiques} $\lambda_1,\ldots,
\lambda_k$ de $f$ sont d\'efinis par (cf \cite{RS}, \cite{DS3}) :
\[ \lambda_l := \lim_{n \to \infty} \left(   \int_X  {f^n}^* \om^l  
\wedge \om^{k-l} \right)^{1/n} . \]

Sous l'hypoth\`ese $d_t >  \lambda_{k-1}$, la suite de mesure 
$\mu_n :=  {1 \over d_t ^n} {f^n}^* \om^k$ converge vers une mesure 
de probabilit\'e invariante ergodique $\mu$ v\'erifiant (cf \cite{RS}, \cite{G}, \cite{DS2}) : \\

($\al$)  Les fonctions \emph{quasi-psh} sur $X$ sont $\mu$-int\'egrables. 

($\beta$)  $\mu$ est de jacobien constant $d_t$ : pour tout $A$ bor\'elien de 
$X\setminus \II$ sur lequel $f$ est injective, on a $\mu(f(A)) = d_t .\mu(A)$. \\

La propri\'et\'e ($\al$) montre que la fonction $\log (\dist (x,\JJ))$ est $\mu$-int\'egrable. 
En effet, elle est minor\'ee par une fonction \emph{quasi-psh} sur $X$. 
C'est clair si $X$ est une vari\'et\'e projective. Si $X$ est seulement 
k\"ahl\'erienne, 
ceci r\'esulte d'un r\'esultat de Blanchard (cf \cite{DS2}, App. A.1). 
Il s'ensuit $\mu(\JJ) = 0$. D'apr\`es le lemme \ref{cauchy} appliqu\'e \`a $f$ et 
$f^{-1}$, 
on peut majorer les valeurs absolues des fonctions $\log \bignorm{ df(x)^{\pm 1} }$ 
et $\log (\Jac f(x))$ par une fonction du type $c_1 - c_2\log (\dist (x,\JJ))$, 
o\`u $c_1$ et $c_2$ sont des constantes positives. 
Par exemple, les deux in\'egalit\'es dans le lemme \ref{cauchy} appliqu\'ees \`a  
$f^{-1}$ donnent
$$\bignorm {df(x)^{-1}} \leq \bignorm {(df^{-1})(f(x))}
\lesssim \dist (f(x), \JJ')^{-p}\lesssim \dist(x,\JJ)^{-q}.$$
Les fonctions 
$\log \bignorm{ df(x)^{\pm 1} }$ et $\log (\Jac f(x))$ sont donc  
$\mu$-int\'egrables et les exposants de Lyapounov  
$\chi_1 \leq \cdots \leq \chi_k$ de $\mu$ sont bien d\'efinis. On dispose de l'\'egalit\'e  :
\[ 2\Sigma := 2(\chi_1 + \cdots + \chi_k)  = \int_X \log \Jac (f) \, d\mu   \]
et de l'in\'egalit\'e $2\Sigma \geq \log d_t$ (cf par exemple \cite{DS1}). Rappelons que les exposants $\chi_1$ et $\chi_k$ sont d\'efinis par : 
\[ - \chi_1 = \lim_{N \to +\infty}  {1 \over N} \int_X \log \bignorm{  df^N(x)  ^{-1} } d\mu(x) = \inf_{N \geq 1} \Big \{ {1 \over N} \int_X \log \bignorm{  df^N(x)  ^{-1} } d\mu(x) \Big\}  \]
et
\[ \chi_k  = \lim_{N \to +\infty}  {1 \over N} \int_X \log \bignorm{  df^N(x) } d\mu(x) = \inf_{N \geq 1} \Big \{ {1 \over N} \int_X \log \bignorm{ df^N(x) } d\mu(x)   \Big\} . \]
Ces \'egalit\'es nous permettent de remplacer $f$ par une it\'er\'ee dans la d\'emonstration du th\'eor\`eme \ref{result}. En effet, le degr\'e topologique de $f^N$ est $d_t^N$ et $\mu$ est aussi la mesure d'\'equilibre de $f^N$, d'exposants $N\chi_1,  \ldots , N\chi_k$. On peut donc supposer que pour $0 < \epsilon_0 < 1$ fix\'e, on a : 
\begin{equation}\label{serd1}
  - \chi_1  \leq \int_X \log \bignorm{  df(x)  ^{-1} } d\mu(x) \leq - \chi_1(1 - \epsilon_0) 
\end{equation}
et 
\begin{equation}\label{serd2}
 \chi_k \leq \int_X \log \bignorm{ df(x) } d\mu(x) \leq \chi_k (1+ \epsilon_0).
\end{equation}


\subsection{Extension naturelle}\label{extnat}


Il s'agit ici de rendre inversible le syst\`eme dynamique $(X , f , \mu)$ en consid\'erant son \emph{extension naturelle} (cf \cite{CFS}, Chap.$10$, §$4$). Posons : 
\[ \widehat {X} := \Big \{ \hat x := (x_n)_{n \in \Z} \, , \, f (x_n) = x_{n+1} \Big\} \subset X^\Z  \]
muni de la topologie produit. On note $\pi_0 : \hat x \mapsto x_0$ la projection ``au temps z\'ero'', et $\hat f$ le d\'ecalage \`a gauche sur $\widehat {X}$ de sorte que $f \circ \pi_0 = \pi_0 \circ \hat f$. L'\emph{extension naturelle} de $(X , f , \mu)$ est le syst\`eme dynamique $(\widehat X , \hat f , \hat \mu)$, o\`u $\hat \mu$ est l'unique mesure de probabilit\'e sur $\widehat {X}$, invariante par $\hat f$, v\'erifiant $\hat \mu (\pi_0 ^{-1} (A)) = \mu(A)$ pour tout bor\'elien $A$ de $X$. Cette mesure est ergodique car $\mu$ est elle-m\^eme ergodique. \\

Observons que l'extension naturelle offre un cadre de travail commode pour d\'efinir les branches inverses des it\'er\'es de $f$. Pour tout $\hat x \in  \widehat X$, on note (lorsqu'elle existe) $f^{-n}_{\hat x}$ la branche inverse de $f^n$ d\'efinie au voisinage de $x_0$, v\'erifiant $f^{-n}_{\hat x}(x_0) = x_{-n}$. Cependant, l'int\'er\^et de l'extension naturelle ne s'arr\^ete pas \`a ce formalisme. L'invariance de $\hat \mu$ par le d\'ecalage \emph{et} son inverse sont cruciales. Cela appara\^it dans la preuve du lemme \ref{var2} (le th\'eor\`eme de Birkhoff est utilis\'e avec $\hat f$ et ${\hat f}^{-1}$) et de mani\`ere cach\'ee dans les d\'emonstrations des in\'egalit\'es du th\'eor\`eme \ref{result}. \\


Nous terminons cette section en introduisant le bor\'elien $\hat f$-invariant 
\[ \widehat X^*  := \left\{ \hat x \in \widehat {X} \, , \,  x_n \notin \JJ \, , \,  \forall n \in \Z \right\} \]
o\`u $\JJ$ a \'et\'e introduit \`a la section \ref{nota}. Cet ensemble est de mesure totale, car $\mu(\JJ) = 0$. Signalons qu'on lui supprimera par la suite des ensembles de mesure nulle, sans pour autant changer de notation.


\subsection{Fonctions \`a variation lente}\label{varle}


La preuve du th\'eor\`eme \ref{result} repose sur l'existence des branches inverses le long d'orbites n\'egatives typiques. Il sera aussi n\'ecessaire de pouvoir comparer les familles de branches inverses $\Big(  f^{-n}_{\hat f^p (x)}  \Big)_{n \geq 0}$ lorsque $p$ devient grand. On s'int\'eresse en particulier aux variations du rayon de la boule sur laquelle celles-ci sont d\'efinies. Cela nous am\`ene \`a introduire la d\'efinition de fonction \`a \emph{variation lente}, qui rejoint celle de fonction \emph{temp\'er\'ee} (\cite{KH}, Chap.S. §2) :  
\begin{defn}\label{exple}
Soit $\epsilon > 0$. Une fonction mesurable $S :  \widehat X \to \R^*_+$ est dite $\epsilon$-lente si il existe une fonction mesurable $\tilde S_1 :  \widehat X \to \R^*_+$ v\'erifiant : 
\[ \pp \forall q \geq 1 \ , \ S(\hat f ^q (\hat x)) \geq e^{-q \epsilon} \tilde S_1 (\hat x). \]
$S$ est dite $\epsilon$-rapide si il existe une fonction mesurable $\tilde S_2 :  \widehat X \to \R^*_+$ telle que :
\[ \pp \forall q \geq 1 \ , \ S(\hat f ^q (\hat x)) \leq e^{q\epsilon} \tilde S_2 (\hat x). \] 
\end{defn}

Observons que le supremum et l'infimum d'une famille finie de fonctions $\epsilon$-rapides (resp.  $\epsilon$-lentes) sont encore des fonctions $\epsilon$-rapides (resp. $\epsilon$-lentes). On peut aussi supposer que $\tilde S_1$ (resp. $\tilde S_2$) prend ses valeurs dans $]0,1]$ (resp. $[1,+\infty[$). Le lemme qui suit repose sur le th\'eor\`eme de Birkhoff.

\begin{lem}\label{var2}
Soient  $\epsilon >0$ et $u : X \to \R^*_+$ une fonction mesurable telle que $\log u \in L^1(\mu)$. Posons $\chi := \int_X \log u \, d\mu$. 
\begin{enumerate}
\item[{\rm (a)}]
Il existe une fonction  $\epsilon$-lente $V_1 : \widehat X \to ]0,1]$, et une fonction $\epsilon$-rapide $V_2 : \widehat X \to [1,+\infty[$ v\'erifiant :
\[ \pp \forall n \geq 1 \ , \ V_1(\hat x)  e^{n(\chi - \epsilon)} \leq \prod_{j=1}^n u(x_{-j}) \leq V_2(\hat x) e^{n(\chi + \epsilon)}. \]  
\item[{\rm (b)}]
Il existe une fonction $\epsilon$-lente $V : \widehat X \to ]0,1]$ v\'erifiant :
\[ \pp \forall n \geq 0 \ , \ u(x_{-n}) \geq V(\hat x) e^{-n\epsilon}. \]  
\end{enumerate}
\end{lem}

\proof
Montrons la premi\`ere partie. Il suffit d'\'etablir l'existence de $V_2$, car celle de $V_1$ s'en d\'eduit en consid\'erant la fonction $1/u$. Le th\'eor\`eme de Birkhoff, appliqu\'e aux syst\`emes $(\widehat X , \hat f , \hat \mu)$, $(\widehat X , {\hat f}^{-1} , \hat \mu)$ et \`a la fonction $\hat x \mapsto \log u (x_0)$, fournit :
\[ \pp \lim_{n \to +\infty} {1 \over n} \sum_{j=1}^{n} \log u(x_{-j}) 
=  \lim_{n \to +\infty} {1 \over n} \sum_{j=0}^{n-1} \log u(x_j) = \chi. \]
Il existe donc une fonction mesurable 
$p :  \widehat X  \to [1,+\infty[$ v\'erifiant :
\[  \pp \forall n \geq 1  \ , \  \prod_{j = 1}^n u(x_{-j}) 
\leq p(\hat x) e^{n(\chi + \epsilon / 2)}, \]
\begin{equation}\label{deux}
 \pp \forall n \geq 0 \ , \ p(\hat x)^{-1} e^{n(\chi - \epsilon / 2)} 
\leq \prod_{j = 0}^n u(x_j) \leq p(\hat x)  e^{n(\chi + \epsilon / 2)}.
\end{equation}
Posons 
\[ V_2'(\hat x) := \sup_{n \geq 1} \Big \{ \prod_{j = 1}^n u(x_{-j})  e^{-n(\chi + \epsilon / 2)} \Big \} \ \textrm{ et } \  V_2(\hat x) := \max \Big \{ 1  , V_2'(\hat x) \Big \} . \]
Ces fonctions sont bien finies presque partout. La fonction $V_2$ v\'erifie l'in\'egalit\'e voulue. Nous montrons maintenant qu'elle est $\epsilon$-rapide. Fixons $\hat x$ 
g\'en\'erique, et posons $p := p(\hat x)$, $V_2 :=  V_2(\hat x)$,  
$\tau_\pm := e^{\chi \pm \epsilon / 2}$ et $\pi_n := \prod_{j = 0}^n u(x_j)$. 
La ligne (\ref{deux}) 
s'\'ecrit alors $p^{-1} \tau_-^n \leq  \pi_n \leq p \tau_+^n$. 
Ainsi, pour tout $q \geq 1$ :
\begin{eqnarray*}
 V_2'(\hat f ^q (\hat x)) &  \leq  &  
\sup \left\{ \frac{\pi_{q-1}\tau_+^{-1}}{\pi_{q-2}} , 
\ldots ,  \frac{\pi_{q-1}\tau_+^{1-q}}{\pi_0 } , 
\pi_{q-1}\tau_+^{-q}  ,  \pi_{q-1} \tau_+^{-q}   V_2 \right\} \\ 
                         & \leq & \sup  \left\{ 
\frac{p^2  \tau_+^{q-2}}{\tau_-^{q-2}} , \ldots ,  
p^2  ,  p  \tau_+^{-1} ,   p  \tau_+^{-1}  V_2  \right \} \\
                         & \leq & \sup  \left\{  p^2  e^{(q-2)\epsilon} , 
\ldots ,  p^2 ,  p \tau_+^{-1} ,  p \tau_+^{-1}  V_2   \right \} \\ 
                         & \leq & p^2  e^{q\epsilon} \max(1,\tau_+^{-1}) V_2 
\end{eqnarray*}                                      
car $p \geq 1$ et $V_2  \geq 1$. Ceci entra\^ine que $V_2$ est $\epsilon$-rapide. 

Le point (b) est une cons\'equence du point (a) en rempla\c{c}ant $\epsilon$ par $\epsilon/2$ : l'in\'egalit\'e est satisfaite avec la fonction $V(\hat x) := \min \{ 1, u(x_0) , e^{\chi} V_1 /V_2 \}$.\fin


\subsection{Dimensions}\label{dimhaus}

   
Nous rappelons la d\'efinition de la dimension de Hausdorff et celle de la dimension d'une mesure. On consultera par exemple les livres de Falconer (\cite{F}, Chap.2) et de Pesin (\cite{P}, Chap.2). Soit $\Lambda$ une partie de $X$. Pour tout $\al \geq 0$ et $\delta > 0$, on note :
\[ \Lambda_{\al,\delta}(\Lambda) :=  \inf_{\GG_\delta} \Big\{ \sum_{U \in \GG_\delta}  (\diam \, U) ^\al \Big \}   \] 
o\`u $\GG_\delta$ parcourt l'ensemble des recouvrements de $\Lambda$ par une famille finie ou d\'enombrable d'ouverts de diam\`etre inf\'erieur \`a $\delta$. La \emph{$\al$-mesure de Hausdorff} et la \emph{dimension de Hausdorff} de $\Lambda$ sont alors d\'efinies par :
\[ \Lambda_\al (\Lambda) := \lim_{\delta \to 0} \uparrow \Lambda_{\al,\delta}(\Lambda), \] 
\[ \dim_\HH (\Lambda) := \inf \Big\{ \al \geq 0 \, , \,  \Lambda_\al (\Lambda) = 0 \Big\} = \sup \Big\{ \al \geq 0 \, , \,  \Lambda_\al (\Lambda) = + \infty \Big\} \]
et la \emph{dimension} de $\mu$ est \'egale \`a :  
\[ \dim_\HH(\mu) :=  \inf \Big\{ \dim_\HH (Z)\, , \, Z 
\textrm{ bor\'elien de } X \, , \, \mu(Z) = 1 \Big\}. \]
Celle-ci nous renseigne sur la mani\`ere dont $\mu$ remplit son support. Observons pour terminer le fait \'el\'ementaire suivant. Nous l'utiliserons pour \'etablir la majoration du th\'eor\`eme \ref{result}.  
\begin{lem}\label{vbnn}
Pour tout bor\'elien $\Lambda$ v\'erifiant $\mu(\Lambda) >0$, on a $\dim_\HH(\mu) \leq \dim_\HH(\Lambda)$.
\end{lem}

\proof
Introduisons l'ensemble 
\[  Z := \Big\{ z \in X \ , \exists n \in \Z \ , \ f^n(z) \cap \JJ \neq \emptyset \Big\} \] 
o\`u $\JJ$ a \'et\'e d\'efini \`a la section \ref{nota}. Puisque $\mu(\JJ) =0$, l'invariance de $\mu$ entra\^ine $\mu(Z) = 0$. L'ensemble invariant $\cup_{n \in \Z} f^n(\Lambda \setminus Z)$ est donc de mesure $1$, car $\mu$ est ergodique. On en d\'eduit :
\[  \dim_\HH (\mu) \leq \dim_\HH \Big( \bigcup_{n \in \Z} f^n(\Lambda \setminus Z) \Big) = \sup_{n \in \Z} \, \dim_\HH f^n(\Lambda \setminus Z)  =  \dim_\HH  (\Lambda \setminus Z) \leq  \dim_\HH  (\Lambda)   \]
o\`u la premi\`ere \'egalit\'e est classique, tandis que la deuxi\`eme provient du caract\`ere m\'eromorphe de $f$. \fin


\section{Branches inverses}\label{brinv}


L'objectif de cette partie est d'\'etablir la proposition \ref{lempr} qui fournit l'existence et le contr\^ole des branches inverses des it\'er\'es de $f$ le long d'orbites g\'en\'eriques. On trouve des r\'esultats similaires dans les articles de Briend-Duval \cite{BD} et Binder-DeMarco \cite{BdM}. Nous traitons ici le cas m\'eromorphe et nous pr\'ecisons comment varient les fonctions $C,\kappa$ et $r$. On se fixe une fois pour toutes $0 < \epsilon \ll \chi_1$. Le r\'eel $\epsilon_0$ et l'ensemble $\widehat X^*$ ont \'et\'e respectivement introduits aux sections \ref{nota} et  \ref{extnat}. On s'autorise aussi \`a supprimer \`a $\widehat X^*$ des ensembles de $\hat \mu$-mesure nulle sans changer de notation. 

\begin{prop} \label{lempr}
Il existe une constante $\rho \geq 1$ ind\'ependante de $\epsilon$ et des fonctions mesurables $C,\kappa : \widehat X^* \to [1,+\infty[$, $r : \widehat X^* \to ]0,1]$ satisfaisant pour tout $\hat x \in \widehat X^*$ :
\begin{enumerate} 
\item[{\rm (a)}]
$r$ est $\rho\epsilon$-lente, $C$ et $\kappa$ sont $\epsilon$-rapides.

\item[{\rm (b)}]
$f^{-n}_{\hat x}$ est d\'efinie sur $B_{x_0} (r (\hat x))$ pour tout $n\geq 1$.

\item[{\rm (c)}]
Pour tout $0 < s \leq r(\hat x)$, on a $B_{x_{-n}} \Big( {s \over C(\hat x)} e^{-n (\chi_k(1+\epsilon_0) + \epsilon) }  \Big) \subset  f^{-n}_{\hat x} \, B_{x_0} (s)$.

\item[{\rm (d)}]Pour tout $0 < s \leq r(\hat x)$, on a 
$\vol \,  f^{-n}_{\hat x} ( B_{x_0} (s))  \leq  \kappa(\hat x) 
\vol B_{x_0}(s). e^{-n ( 2 \Sigma - \epsilon  ) }$.

\end{enumerate}
\end{prop}

La d\'emonstration repose sur les deux lemmes suivants. Le premier montre que les orbites n\'egatives typiques de $f$ ne s'approchent pas trop vite de l'ensemble analytique $\JJ$ d\'efini \`a la section \ref{nota} :

\begin{lem}\label{dfc}
Il existe une fonction $\epsilon$-lente $A_1 : \widehat X^* \to ]0,1]$ telle que :
\[ \forall \, \hat x \in  \widehat X^* \ , \ \forall n \in \N \ , \ \dist \left( x_{-n} , \JJ \right)  \geq 2A_1 (\hat x) e^{-n\epsilon}. \] 
\end{lem}

\proof
La fonction $u(x) := \log \dist ( x , \JJ )$ est $\mu$-int\'egrable, d'apr\`es la section \ref{nota}. On invoque alors le lemme  \ref{var2}-(b). On supprime ici \`a l'ensemble $\widehat X^*$ un ensemble de mesure $\widehat\mu$ nulle.  \fin

Le deuxi\`eme lemme est un r\'esultat d'inversion locale quantitative. Pour tout $\hat x \in \widehat X^*$, on pose : 
\begin{equation*}
 \VV_n(\hat x) := \left \{ z \in X \, , \, \dist(z,\JJ) \geq  A_1 (\hat x) e^{-n\epsilon} \right\} .
\end{equation*}
On introduit \'egalement (cf section \ref{nota}) :
\begin{equation*}
M_n(\hat x) :=  1 +   \bignorm {d f} _{{\mathcal V} _{n+1}(\hat x)} +   \bignorm {d^2 f} _{{\mathcal V} _{n+1}(\hat x)}  +  \bignorm {d f^{-1}} _{{\mathcal V} _{n+1}(\hat x)} +  \bignorm {d^2 f^{-1}} _{{\mathcal V} _{n+1}(\hat x)}. 
\end{equation*}
Le lemme \ref{cauchy} appliqu\'e \`a $f$ et $f^{-1}$ avec $\AA := \II$ ou $\JJ'$ fournit des constantes $c,c' \geq 1$ et $\rho' \geq  1$ ind\'ependantes de $\epsilon$ v\'erifiant :
\begin{equation*} \label{no2}
 \forall \, \hat x \in  \widehat X^*  \ , \  \exists \, z \in  {\mathcal V}_{n+1}(\hat x)  \ , \  M_n(\hat x)^{2k+1} \leq c' \textrm{dist} (z,{\mathcal J})^{-\rho'}  \leq   c\left( A_1 (\hat x) e^{-(n+1)\epsilon} \right)^{-\rho'}.  
\end{equation*}
Pour tout $\hat x \in  \widehat X^*$ et tout $n \in \N$, notons :
\[ \al_n(\hat x) := \min   \left \{   1 \, , \, \bignorm{ d f (x_{-n-1})  }  \,  , \,  \bignorm{ d f (x_{-n-1}) ^{-1}  } \, , \,  \Jac f (x_{-n-1}) \right \} .  \]
On pose alors : 
\[ r_n (\hat x) : = {  \alpha_n(\hat x)  A_1 (\hat x) e^{-(n+1)\epsilon} \over C M_n(\hat x)^{2k+1}  }   \]
o\`u $C$ est une constante v\'erifiant  $C > \max \{ e^{\epsilon/2} , ( 1 - e^{-\epsilon/2} )^{-1} \}$. Le th\'eor\`eme d'inversion locale ``pr\'ecis\'e'' est le suivant.

\begin{lem}\label{bbb} Soit  $B_n := B_{x_{-n}}(r_n (\hat x))$. Alors, pour tout $n \geq 0$, 
\begin{enumerate}
\item[{\rm (a)}]
$f$ poss\`ede une branche inverse $g$ sur $B_n$ v\'erifiant $g(x_{-n}) = x_{-n-1}$.

\item[{\rm (b)}]
$\Lip ( g _ {\vert B_n}) \leq \bignorm{ df(x_{-n-1}) ^{-1} } e ^{\epsilon/2}$.

\item[{\rm (c)}]
$\Lip ( f _ {\vert g (B_n)}) \leq \bignorm{ df(x_{-n-1}) }e ^{\epsilon/2}$.

\item[{\rm (d)}]   
$\inf  \Big\{  \Jac f (y)  , y \in g(B_n)  \Big\} \geq  \Jac f (x_{-n-1}) e ^{-\epsilon/2}$.

\item[{\rm (e)}] 
Il existe une fonction $\eta : \widehat X^* \to ]0,1]$, qui est $3\rho'\epsilon$-lente, et telle que 
\[ \forall n \geq 0 \ , \ r_n(\hat x) \geq e^{-3 n \rho' \epsilon} \eta(\hat x). \]
\end{enumerate}
\end{lem}

\proof

\noindent (a) On a $B_n \subset \VV_n(\hat x) \subset X \setminus \JJ$ car $r_n(\hat x) \leq  A_1 (\hat x) e^{-n\epsilon} \leq d(x_{-n},\JJ)/2$ (cf lemme \ref{dfc}). On termine en rappelant que $\JJ$ contient les valeurs critiques de $f$.

\noindent (b) Par d\'efinition de $M_n(\hat x)$ et de $r_n(\hat x)$, on obtient en utilisant l'in\'egalit\'e des accroissements finis :
\begin{eqnarray*}
   \Lip ( g _ {\vert B_n})      & \leq  & \sup  \, \{  \bignorm{ dg(y) } \, , \, y \in B_n  \}   \\ 
                                & \leq  & \bignorm{ df(x_{-n-1}) ^{-1} } + r_n(\hat x) M_n(\hat x)    \\
                                & \leq  &  \bignorm{ df(x_{-n-1}) ^{-1} } +   \bignorm{ df(x_{-n-1}) ^{-1} } / C   \\ 
                                & \leq  &  \bignorm{ df(x_{-n-1}) ^{-1} } e^{\epsilon /2}
\end{eqnarray*}
car $C^{-1} \leq  1 - e^{-\epsilon/2} \leq   e^{\epsilon/2} -1$.

\noindent (c) D'apr\`es (b) , $g(B_n)$ est contenu dans la boule centr\'ee en $x_{-n-1}$ et de rayon :
\[ r_n'(\hat x) := r_n(\hat x)   \bignorm{ df(x_{-n-1}) ^{-1} } e^{\epsilon /2}  \leq  r_n(\hat x) M_n(\hat x) e^{\epsilon /2}. \]
Elle est en particulier contenue dans $\VV_{n+1}(\hat x)$. En effet, on a :
\[  r_n(\hat x) M_n(\hat x) e^{\epsilon /2} \leq A_1(\hat x) e^{-(n+1)\epsilon}  e^{\epsilon /2} /C \leq  A_1(\hat x) e^{-(n+1)\epsilon}  \leq  d(x_{-n-1},{\mathcal J})/2.  \]
On en d\'eduit comme pr\'ec\'edemment :
\begin{eqnarray*} 
 \textrm{Lip} ( f _ {\vert g(B_n)})       &  \leq  &   \bignorm{ df(x_{-n-1}) }  +   r_n'(\hat x) M_n(\hat x) \\
                                 &  \leq  &   \bignorm{ df(x_{-n-1}) }  +   r_n(\hat x) M_n^2(\hat x) e^{\epsilon /2}  \\
                                 &  \leq  &   \bignorm{ df(x_{-n-1}) }  +  \bignorm{ df(x_{-n-1}) } / C  \\
                                 &  \leq  &   \bignorm{ df(x_{-n-1}) } e^{\epsilon /2} .
\end{eqnarray*}

\noindent (d) On a dans ce cas :
\noindent (d) On a dans ce cas :
\begin{eqnarray*} 
    \inf  \Big\{  \textrm{Jac} f (y)  , y \in g(B_n)  \Big\}    &  \geq  &    \textrm{Jac} f (x_{-n-1}) -  r_n'(\hat x) M_n(\hat x)^{2k} \\
                                                                &  \geq  &    \textrm{Jac} f (x_{-n-1}) -  r_n(\hat x) M_n(\hat x)^{2k+1} e^{\epsilon /2}  \\
                                                               &  \geq  &    \textrm{Jac} f (x_{-n-1}) -  \textrm{Jac} f (x_{-n-1}) / C  \\
&  \geq  &   \textrm{Jac} f (x_{-n-1}) e^{-\epsilon /2}.
\end{eqnarray*}

\noindent (e) Puisque les fonctions \emph{quasi-psh} sont $\mu$-int\'egrables, il en est de m\^eme pour 
\[  z \mapsto  \log \left(  \min   \left \{   1 \, , \, \bignorm{ d f (z)  }  \,  , \, \bignorm{ d f (z) ^{-1}  } \, , \,  \Jac f (z) \right \} \right) . \] 
Ainsi, d'apr\`es le lemme \ref{var2}(b), il existe une fonction $\epsilon$-lente $A_2 : \widehat X \to ]0,1]$ v\'erifiant 
\[ \pp \forall n \geq 0 \ , \ \al_n(\hat x)  \geq A_2(\hat x) e^{-n\epsilon} .\] 
A l'aide de la ligne (\ref{no2}), on v\'erifie facilement que l'on a (avec $\rho' \geq 1$) :
\[ \pp \forall n \geq 0 \ , \ r_n(\hat x)  \geq e^{-3 n \rho' \epsilon}  \left( (Cc)^{-1}  e^{-\epsilon}  A_1(\hat x)^{2\rho'} A_2(\hat x)^{\rho'} \right).   \] 
On prend alors pour $\eta$ la fonction entre les parenth\`eses : elle est bien $3\rho'\epsilon$-lente car $A_1$ et $A_2$ sont $\epsilon$-lentes.  \fin

\noindent \textsc{D\'emonstration de la proposition} \ref{lempr} : \rm  
Les notations sont celles du lemme \ref{bbb}. Posons 
\[ \forall n \geq 1 \ , \  \tilde B_n := \left \{ f^{n-1}(z) \, , \, z \in B_{n-1} \, ,  \, f(z) \in B_{n-2} \,  ,  \, \ldots \,  ,  \, f^{n-1}(z) \in B_0    \right\} \subset B_0 .   \]
D'apr\`es le lemme \ref{var2}, le lemme \ref{bbb}, et les in\'egalit\'es des lignes (\ref{serd1}) et (\ref{serd2}) (cf section \ref{nota}), il existe des fonctions $\epsilon$-rapides $C_1, C_2 :  \widehat X^* \to [1,+\infty[$ et une fonction $\epsilon$-lente $C_3 : \widehat X^* \to ]0,1]$ v\'erifiant : 
\begin{enumerate}
\item[(1)] 
Pour tout $n \geq 1$, $g_n := f^{-n}_{\hat x}$ existe sur $\tilde B_n$.

\item[(2)]
$\Lip (  g_{n \, \vert \tilde B_n } ) \leq  \prod_{i=1}^n \bignorm{ df(x_{-i}) ^{-1} } e ^{n \epsilon/2} \leq  C_1(\hat x) e ^{n (-\chi_1(1-\epsilon_0) + \epsilon)} $.

\item[(3)]
$\Lip \Big( f^n_{\vert g_n (\tilde B_n)} \Big) \leq \prod_{i=1}^n \bignorm{ df(x_{-i}) }e ^{n\epsilon/2} \leq  C_2(\hat x) e ^{n(\chi_k(1+\epsilon_0) + \epsilon)} $.

\item[(4)]   
$\inf \left\{ \Jac f^n (y) , y \in  g_n(\tilde B_n) \right\} \geq \Jac f^n (x_{-n}) e ^{-n \epsilon/2} \geq C_3(\hat x) e^{n(2\Sigma - \epsilon)}$.
\end{enumerate}

Pour le dernier point, rappelons que l'on a $\Jac f^n (x_{-n}) = \prod_{i=1}^n \Jac f (x_{-i})$ et $\int_X \log \Jac (f) d\mu = 2 \Sigma$ (cf section \ref{nota}). Posons $C := \max \{ C_1,C_2 \}$ et $r := \eta / C$, o\`u $\eta$ a \'et\'e introduite au lemme \ref{bbb}. Ces fonctions sont respectivement $\epsilon$-rapide et  $\rho \epsilon$-lente, avec $\rho := 3\rho'+1$. Pour montrer (b), on v\'erifie l'inclusion $B(x_0,r(\hat x)) \subset \tilde B_n$ par r\'ecurrence, en utilisant (2) avec $\epsilon \ll \chi_1$. L'assertion (c) d\'ecoule de (3). Enfin, le point (d) provient du th\'eor\`eme de changement de variables avec $\kappa := C_3^{-1}$. \fin


\section{D\'emonstration du th\'eor\`eme \ref{result}}\label{demo}


\subsection{Minoration} 

\'Etant donn\'e $\epsilon >0$, prenons les fonctions \`a variation lente $r,C$ d\'efinies par la proposition \ref{lempr}. La fonction $\zeta := r / C$ est donc $(\rho +1)\epsilon$-lente, et il existe une fonction $\si : \widehat X^* \to ]0,1]$ v\'erifiant :
\[ \forall n \geq 1 \ , \ \zeta(\hat f ^n (\hat x)) \geq \si(\hat x).e^{-n(\rho+1)  \epsilon}. \]
On obtient ainsi, en posant $\delta_n := e^{-n(\chi_k(1+\epsilon_0) + \epsilon)}.e^{-n(\rho+1)\epsilon}$ : 
\[ \pp  \forall n \geq 1 \ , \ B_{x_0} (\si(\hat x).\delta_n) \subset B_{x_0} ( \zeta(\hat f^n (\hat x)) . e^{-n(\chi_k(1+\epsilon_0) + \epsilon)}) .  \]
Puisque $f^n$ est injective sur ces boules (cf proposition \ref{lempr}-(c)), et que $\mu$ est une mesure de \emph{probabilit\'e} et de \emph{jacobien} $d_t$ (cf section \ref{nota}), cette inclusion implique :
\[ \pp  \forall n \geq 1 \ , \ \mu \left(  B_{x_0} (\sigma (\hat x).\delta_n) \right) \leq d_t^{-n}  \]
ce qui fournit imm\'ediatement :
\[ \pp  \forall n \geq 1 \ , \ { \log \mu \left(  B_{x_0} (\sigma (\hat x).\delta_n) \right)  \over  \log ( \sigma (\hat x).\delta_n )   }    \geq { \log  d_t^{-n} \over \log ( \sigma (\hat x).\delta_n ) }.  \]
On en d\'eduit, en prenant la limite inf\'erieure :
\[ \pp \liminf_{n \to +\infty} { \log \mu \left(  B_{x_0} (\sigma (\hat x).\delta_n) \right)  \over  \log ( \sigma (\hat x).\delta_n )   }    \geq { \log  d_t \over \chi_k(1+\epsilon_0) + (\rho+2)\epsilon }.  \] 
Il existe donc un bor\'elien $\Lambda$ de $X$ v\'erifiant $\mu(\Lambda) = 1$, et $\theta : \Lambda \to ]0,1]$ tels que : 
\begin{equation}\label{eee}
 \forall x \in \Lambda \ , \ \liminf_{n \to +\infty} { \log \mu \left(  B_x (\theta (x).\delta_n) \right)  \over  \log ( \theta (x).\delta_n )   }    \geq { \log  d_t \over \chi_k(1+\epsilon_0) + (\rho+2)\epsilon }. 
\end{equation} 
Le corollaire \ref{bore} appliqu\'e \`a $\rho_n(x) := \theta(x).\delta_n$ entra\^ine :
\[ \dim_\HH(\mu) \geq { \log d_t \over \chi_k(1+\epsilon_0) +  (\rho+2) \epsilon }. \]
On obtient la minoration $\dim_\HH(\mu) \geq \log d_t / \chi_k$ en faisant tendre $\epsilon$ et $\epsilon_0$ vers z\'ero.

\subsection{Majoration}

Consid\'erons pour $\epsilon >0$ fix\'e les fonctions \`a variation lente $r,C$ et $\kappa$ d\'efinies \`a la proposition \ref{lempr}. On leur associe des fonctions mesurables $\tilde r$, $\tilde C$ et $\tilde \kappa$ v\'erifiant les propri\'et\'es de la d\'efinition \ref{exple}. Fixons aussi $r_0,C_0,\kappa_0 > 0$ pour que le bor\'elien   
\[ \widehat A := \Big \{ \hat x \in \widehat X^* \, , \, \tilde r(\hat x) \geq r_0 \, , \, \tilde C(\hat x) \leq C_0 \, , \, \tilde \kappa(\hat x) \leq \kappa_0 \Big \} \]
v\'erifie $\hat \mu(\widehat A) >0$. On note alors $A := \pi_0 (\widehat A)$, et pour tout $n \geq 0$ :
\[ r_n := r_0 e^{-n \rho \epsilon} \ , \ C_n := C_0 e^{n \epsilon} \ , \  \kappa_n := \kappa_0 e^{n \epsilon}  \ , \ \delta_n :=  (r_n/C_n). e^{-n ( \chi_k(1+\epsilon_0) + \epsilon  )}.  \]
Pour tout $\hat x \in \widehat A$ et $n \geq 0$, on pose aussi : 
\[ \hat x_n := \hat f ^n (\hat x)  \ , \    g_n := f^{-n}_{\hat x_n}  \ , \ x := \pi_0(\hat x)  \ , \   x_n := \pi_0 (\hat x_n)=f^n(x).   \]  
Le lemme suivant est une cons\'equence de la proposition \ref{lempr} :
\begin{lem}\label{enplus}
Pour tout $\hat x \in \widehat A$, on a : 
\begin{enumerate}
\item[{\rm (a)}]
$g_n$ est d\'efinie sur $B_{x_n}(r_n)$.

\item[{\rm (b)}]
Pour tout $s \in ]0, r_n]$, $B_x ({s \over r_n} \delta_n ) \subset g_n(B_{x_n}(s))$.

\item[{\rm (c)}]
$\vol g_n(B_{x_n}(r_n)) \leq  c e^{-2n\Sigma}$, o\`u $c$ est une constante ind\'ependante de $n$.
\end{enumerate}
\end{lem}

\proof
Les deux premiers points s'obtiennent \`a l'aide de la proposition \ref{lempr}-(b),(c), en observant :
\[ r(\hat x_n) \geq \tilde r(\hat x) e^{-n \rho \epsilon} \geq r_0  e^{-n \rho \epsilon} = r_n \ \textrm{ et } \ C(\hat x_n) \leq \tilde C(\hat x) e^{n \epsilon} \leq C_0  e^{n \epsilon} = C_n . \]
Pour le dernier point, la  proposition \ref{lempr}-(d) et l'in\'egalit\'e $\kappa (\hat x_n) \leq \kappa_0  e^{n \epsilon}$ entra\^inent :
\begin{equation*}
 \vol g_n(B_{x_n}(r_n))  \leq  \kappa_0 e^{n \epsilon} \vol (B_{x_n}(r_n)) e^{-n ( 2 \Sigma - \epsilon) } \leq c e^{-2n\Sigma} 
\end{equation*}
avec $c$ une constante assez grande. On utilise ici l'in\'egalit\'e $\textrm{vol} (B_{x_n}(r_n)) \leq e^{-2n \epsilon}$ valable \`a une constante multiplicative pr\`es. \fin
 
Nous obtenons la majoration en montrant l'in\'egalit\'e :
\begin{equation} \label{1}
 \dim_\HH (A)  \, \leq   \, \alpha_\epsilon \, := \, 2k - { 2 \Sigma - \log d_t -2k\rho\epsilon  \over \chi_k(1+\epsilon_0) + (\rho+2)\epsilon}.   \end{equation}
En effet, si celle-ci est v\'erifi\'ee, le lemme \ref{vbnn} entra\^ine $\dim_\HH (\mu) \leq \al_\epsilon$ et on termine en faisant tendre $\epsilon_0$ et $\epsilon$ vers z\'ero. On s'attache \`a pr\'esent \`a la preuve de (\ref{1}). \'Etant donn\'e $\delta >0$, fixons $n \geq 0$ pour que $2\delta_n < \delta$ et choisissons $E_n$, un sous-ensemble fini de $A$, v\'erifiant :

(i) les boules de la famille $\left \{ B_x (\delta_n / 4) , x \in E_n \right \}$ sont disjointes deux \`a deux.

(ii) $A$ est recouvert par la r\'eunion $\bigcup_{x \in E_n}  B_x(\delta_n)$.

\noindent On a alors par d\'efinition de la fonction d'ensemble $\Lambda_{\al,\delta}$ (cf section \ref{dimhaus}) :
\begin{equation}\label{aa}
  \Lambda_{\al_\epsilon,\delta}(A) \leq  \# E_n . (2 \delta_n)^{\al_\epsilon} . 
\end{equation}
Nous montrons maintenant l'estimation (\ref{1}) en v\'erifiant que le membre de droite de (\ref{aa}) est born\'e par une constante ind\'ependante de $\delta$. Nous \'evaluons $\# E_n$ \`a l'aide du r\'esultat suivant :
\begin{lem}\label{truc}
Il existe une famille d'ouverts $\PP_n$ de cardinal au plus $\si r_n^{-2k} d_t^n$, o\`u $\si$ est une constante ind\'ependante de $n$, v\'erifiant :
\begin{enumerate}
\item[{\rm (a)}]
$\forall P \in \PP_n \, , \, \vol (P) \leq c e^{-2n\Sigma}$.

\item[{\rm (b)}]
$\forall x \in E_n \, , \,  B_x(\delta_n / 4) \subset \bigcup_{P \in \PP_n} P$.
\end{enumerate}
\end{lem}
Supposons ce r\'esultat acquis, et terminons la d\'emonstration du th\'eor\`eme. Si $\Ga_n$ d\'esigne le volume d'une boule de rayon $\delta_n / 4$, le lemme \ref{truc} et (i) entra\^inent :
\begin{equation}\label{bb}
\# E_n  \leq   {1 \over \Ga_n }  \sum_{P \in \PP_n} \vol(P) \leq  {1 \over \Ga_n } \si r_n^{-2k} d_t^n . c  e^{-2n\Sigma}  .
\end{equation}
En utilisant (ii) et les lignes (\ref{aa}) et (\ref{bb}), on obtient les in\'egalit\'es suivantes (\`a des constantes multiplicatives pr\`es ind\'ependantes de $n$) :
\begin{eqnarray*}
   \Lambda_{\al_\epsilon,\delta} (A)   & \lesssim  & \# E_n . \delta_n^{\al_\epsilon} \\ 
                                       & \lesssim  & \delta_n ^{\al_\epsilon - 2k}  r_n^{-2k} e^{-n(2\Sigma - \log {d_t})}   \\
                                       & \lesssim  & \left( {r_n \over C_n} \right) ^{\al_\epsilon - 2k}    e^{-n ( \chi_k(1+\epsilon_0) + \epsilon) (\al_\epsilon - 2k) } r_n^{-2k}  e^{-n(2\Sigma - \log {d_t})}  \\ 
                                       & \lesssim  & e ^{n \left[       \epsilon (\rho + 1) + ( \chi_k(1+\epsilon_0) + \epsilon )  \right]    \left( { 2\Sigma - \log {d_t} - 2k\rho \epsilon   \over  \chi_k(1+\epsilon_0) + (\rho+2)\epsilon } \right)  }  e ^{2kn \rho \epsilon} e^{-n(2\Sigma -  \log {d_t})} = 1 .
\end{eqnarray*}
Ainsi, $\Lambda_{\al_\epsilon,\delta} (A)$ est major\'e par une constante ind\'ependante de $\delta$. La quantit\'e $\Lambda_{\al_\epsilon} (A)$ est donc finie, ce qui entra\^ine (\ref{1}).  \fin 

\noindent \textsc{D\'emonstration du lemme} \ref{truc} : \rm
Soient $n \geq 0$ et $\BB$ une famille finie de boules recouvrant $X$ de rayon $r_n/4$ de cardinal major\'e par $\si r_n^{-2k}$, o\`u $\si$ est une constante ind\'ependante de $n$. Reprenons les notations et les r\'esultats du lemme \ref{enplus}. Soient $x \in A$ et $B \in \BB$ tels que $f^n(x) \in B$. Si $2B$ d\'esigne la boule concentrique \`a $B$ et de rayon double, on a clairement :
\begin{equation}\label{incinc}
 B_{x_n}({r_n / 4})  \subset 2B \subset B_{x_n}({r_n})  .
\end{equation}
Posons $P_x := g_n(2B)$. L'inclusion de droite dans (\ref{incinc}) et le lemme \ref{enplus}-(a) montrent que l'ouvert $P_x$ est bien d\'efini et qu'il est contenu dans $g_n(B_{x_n}(r_n))$. Le lemme \ref{enplus}-(c) entra\^ine : 
\[ \vol (P_x) \leq  c e^{-2n\Sigma}  . \]
Si $x \in E_n$, l'inclusion de gauche dans (\ref{incinc}) et le lemme \ref{enplus}-(b) avec $s = r_n /4$ impliquent : 
\[ B_x(\delta_n / 4 ) \subset g_n(B_{x_n}({r_n / 4}))  \subset P_x  . \]
La famille $\PP_n := \{ g_n(2B) \, , \, B \in \BB \}$ r\'epond au probl\`eme pos\'e. Son cardinal est bien major\'e par $\si  r_n^{-2k} d_t^n$ ($f^n$ est de degr\'e topologique $d_t^n$) et les applications $g_n$ sont des branches inverses de $f^n$ sur les boules de $\BB$. \fin

\vspace{1 cm}

Tien-Cuong Dinh et Christophe Dupont

B\^at. 425, Math\'ematique, UMR 8628

Universit\'e Paris-Sud

91405 Orsay Cedex, France. \\

TienCuong.Dinh@math.u-psud.fr

Christophe.Dupont@math.u-psud.fr

\end{document}